\documentclass[12pt]{amsart}
\usepackage[cp1250]{inputenc}
\usepackage[T1]{fontenc}
\usepackage{amscd}

\def\r{\rightarrow}
\def\t{\times}
\def\s{\subset}

\def\B{\mathcal{B}}
\def\C{\mathcal{C}}

\def\H{\mathcal{H}}

\def\L{\mathcal{L}}

\def\t{\times}
\def\ci{\circ}
\def\rz{\Bbb R}

\def\ld{,\ldots ,}

\def\12{{1\over 2}}
\def\rz{\mathbb{R}}

\DeclareMathOperator{\supp}{supp}
\DeclareMathOperator{\id}{id}

 \theoremstyle{plain}
\newtheorem{thm}{Theorem}[section]
\newtheorem{lem}[thm]{Lemma}

\newtheorem{cor}[thm]{Corollary}

\theoremstyle{definition}

\newtheorem{rem}[thm]{Remark}

\keywords{group of homeomorphisms, equivariant homeomorphism,
perfectness, commutator .  } \subjclass{22A05, 57S05}
\thanks{Supported by the AGH grant n. 11.420.04}
\address{Faculty of Applied Mathematics, AGH University of Science and
\linebreak Technology, al. Mickiewicza 30, 30-059 Krak\'ow,
Poland} \email{tomasz@uci.agh.edu.pl}

\title{On commutators of equivariant homeomorphisms}
\author{ Tomasz Rybicki}

\begin{document}

\maketitle

\begin{abstract}
The long-standing problem of the perfectness of the
 compactly supported equivariant homeomorphism group on a
$G$-manifold (with one orbit type) is solved in the affirmative.
The proof is based on an argument different than that for the case
of diffeomorphisms. The theorem is a starting point for computing
$H_1(\H_G(M))$ for more complicated $G$-manifolds.

\end{abstract}

\section{Introduction}

Let $M$ be a topological manifold and let $\H(M)$ be the identity
component of the group of all compactly supported  homeomorphisms
of $M$ with the compact-open topology. J.N. Mather [8] showed that
the group $\H(\rz^m)$ is acyclic. His result combined with a
fragmentation property for homeomorphisms (Corollary 3.1 in [7])
yields that $\H(M)$ is perfect and it is simple as well (cf.[5]).
Recall that a group $H$ is called perfect if $H=[H,H]$, where
$[H,H]$ is the commutator subgroup of $H$.

 Let $G$ be a
compact  Lie group  acting on $M$. Let $\H_G(M)$ be the group of
all equivariant homeomorphisms of $M$ which are isotopic to the
identity through compactly supported equivariant isotopies.
Suppose now that $G$ acts freely on $M$. Then $M$ can be regarded
as the total space of a principal $G$-bundle $\pi:M\r \bar M=M/G$
(cf. [6]).

\begin{thm} $\H_G(M)$ is a perfect group.
\end{thm}

An analogous theorem for equivariant $C^r$-diffeomorphisms, where
$r=1,\ldots,\infty$, $r\neq \dim(\bar M)+1$, is due to A. Banyaga
([4] or [5])  for $G$ being a torus, and is due to K. Abe and K.
Fukui [1] for an arbitrary compact Lie group $G$. Recently the
latter authors showed in [2] such a theorem for equivariant
Lipschitz homeomorphisms. The theorem and the corollary below seem
to be a starting point for computing $H_1(\H_G(M))$ for more
complicated $G$-manifolds, e.g. analogously as it was done in [3]
for equivariant diffeomorphisms.

\section{A clue lemma}
First we prove a clue auxiliary result.  We define a countable set
of intervals of $ \rz^m$, $I_n^k$, where $n=0,1,\ldots$, and
$k=0,1,\ldots,2^n-1$. We set $I^0_0=(\frac 1 3,\frac 2 3
)\t(0,1)^{m-1}$, $I^0_1=(\frac 7 6,\frac 8 6)\t(0,1)^{m-1}$,
$I^1_1=({10\over 6},{11\over 6})\t(0,1)^{m-1}$, $I^0_2=({25\over
12},{26\over  12})\t(0,1)^{m-1}$, $I^1_2=({28\over 12},{29\over
12})\t(0,1)^{m-1}$, $I^2_2=({31\over 12},{32 \over
12})\t(0,1)^{m-1}$, $I^3_2=({34\over 12},{35\over
12})\t(0,1)^{m-1}$. In general we let
$$
I^k_n=\left(n+{1+3k\over 3\cdot 2^n}, n+{2+3k\over 3\cdot
2^n}\right)\t(0,1)^{m-1}.
$$
Let $h$ be an increasing homeomorphism of $(0,1)$ onto
$(0,\infty)$. We define the family $\B=\{B^k_n\}$ of open balls
contained in $(0,1)^m\subset \rz^m$
$$ B^k_n=(h\t\id_{\rz^{m-1}})^{-1}(I^k_n).
$$
  Next, for  $n=0,1,\ldots$, we can choose a family $\{\phi^-_{n}\} $
 of compactly supported homeomorphisms
with pairwise disjoint supports such that
$$
\phi^-_{n}(B^k_{2n})=B^{k}_{2n+1},\quad k=0,1\ld 2^{2n}-1,
$$
and $ \phi^-_{n}=\id$ on
$\bigcup_{\B-\{B^k_{2n},B^{k}_{2n+1}\}}B^{k'}_{n'}$.
 Likewise, for  $n=0,1,\ldots$, we
can find three next families of compactly supported homeomorphisms
with pairwise disjoint supports: the family $\phi^+_{n}$
satisfying
$$ \phi^+_{n}(B^k_{2n})=B^{2^{2n}+k}_{2n+1},\quad k=0,1\ld 2^{2n}-1,
$$
and $ \phi^+_{n}=\id$ on
$\bigcup_{\B-\{B^k_{2n},B^{2^{2n}+k}_{2n+1}\}}B^{k'}_{n'}$; the
family $\psi^-_{n}$ given by
$$
\psi^-_{n}(B^k_{2n+1})=B^{k}_{2n+2},\quad k=0,1\ld 2^{2n+1}-1,
$$
and $ \psi^-_{n}=\id$ on
$\bigcup_{\B-\{B^k_{2n+1},B^{k}_{2n+2}\}}B^{k'}_{n'}$; and the
family $\psi^+_{n}$ verifying
$$
\psi^+_{n}(B^k_{2n+1})=B^{2^{2n+1}+k}_{2n+2},\quad k=0,1\ld
2^{2n+1}-1,
$$
and $ \phi^+_{n}=\id$ on
$\bigcup_{\B-\{B^k_{2n+1},B^{2^{2n+1}+k}_{2n+2}\}}B^{k'}_{n'}$.

 We choose any $\phi^-_{n}$ with support in a  possibly small
ball containing the union of $B^k_{2n}$ and $B^{k}_{2n+1}$,
$k=0,1\ld 2^{2n}-1$, and similarly for the other families. Put
$\phi^*=\prod_{n} \phi^*_{n}$, and $\psi^*=\prod_{n}\psi^*_{n}$,
where $*$ is equal to $-$ or $+$. It is easily seen that then
$\phi^*$ and $\psi^*$ are homeomorphisms.

Let $\C_c(\rz^m)$ is the space of all compactly supported
$\rz$-valued continuous functions on $\rz^m$, and
$\C_B(\rz^m)\s\C_c(\rz^m)$ be the set of functions compactly
supported in $B=B^0_0$. Suppose $u=u^0_0\in\C_B(\rz^m)$.
 We shall define a family of continuous functions $u_n^k\in \C_c(\rz^m)$, where
$n=0,1,\ldots$, and $k=0,1,\ldots,2^n-1$, such that
$\supp(u_n^k)\s B^k_n$.   Set $u_1^0=-\12
u_0^0\circ(\phi^-_{0})^{-1}$ and $u_1^1=-\12
u_0^0\circ(\phi^+_{0})^{-1}$. Next we put $u_2^0=-\12
u^0_1\circ(\psi^-_{0})^{-1}$, $u_2^1=-\12 u^1_1\circ
(\psi^-_{0})^{-1}$, $u^2_2=-\12 u^0_1\circ(\psi^+_{0})^{-1}$, and
$u_2^3=-\12 u^1_1\circ (\psi^+_{0})^{-1}$. In general we define by
induction
$$
u_{2p+1}^{k}=-\12 u_{2p}^k\circ
(\phi^-_{p})^{-1}\quad\hbox{and}\quad
 u_{2p+1}^{2^{2p+1}+k}=-\12
u_{2p}^k\circ (\phi^+_{p})^{-1},
$$
where $p=0,1,\ldots ;k=0,1,\ldots,2^{2p}-1$, and
$$
u_{2p+2}^{k}=-\12 u_{2p+1}^k\circ
(\psi^-_{p})^{-1}\quad\hbox{and}\quad u_{2p+2}^{2^{2p+1}+k}=-\12
u_{2p+1}^k\circ (\psi^+_{p})^{-1},
$$
where $p=0,1,\ldots ;k=0,1,\ldots,2^{2p+1}-1$.

Next we consider four  continuous functions
$$v_1^-=\sum u_{2p+1}^{k},\quad
v_1^+=\sum u_{2p+1}^{2^{2p+1}+k},$$
 where   $p= 0,1,\ldots$,
$k=0,1\ld2^{2p}-1$, and
$$ v_2^-=\sum
u_{2p+2}^{k},\quad
 v_2^+=\sum u_{2p+2}^{2^{2p+1}+k},$$
 where $ p= 0,1,\ldots$, $k=0,1\ld 2^{2p+1}-1$.  Now if we put $\bar w=\sum_{n\geq 0;k}u_n^k$
and $w=\sum_{n\geq 1;k}u_n^k$ then obviously $u=u^0_0=\bar w-w$.
On the other hand, $ w=v^-_2-v^-_2\circ\psi^- +v^+_2
-v^+_2\circ\psi^+$ and, similarly, $\bar w=v^-_1-v^-_1\circ\phi^-
+v^+_1 -v^+_1\circ\phi^+ $.

We shall apply the above reasoning to the semi-direct product
group $\H(\rz^m)\t_{\tau}\C_c(\rz^m)$, where
 $\tau_h(u)=u\circ h^{-1}$ for $h\in
\H(\rz^m)$ and $u\in \C_c(\rz^m)$. Then we have
$$(h_1,u_1)\cdot(h_2,u_2)=(h_1\circ h_2, u_1\circ h_2^{-1}+u_2)$$
for all $h_1,h_2\in\H(\rz^m)$ and $u_1,u_2\in \C_c(\rz^m)$.

Take $(h,\bar u)\in\H(\rz^m)\t_{\tau}\C_c(\rz^m)$. We have
$(h,\bar u)=(\id,u)\cdot(h,0)$, where $u=\bar u\circ h$. Then
$(h,0)$ belongs to the commutator subgroup by [8]. We may assume
that $u$ is supported in $B$ and $u=\bar w-w$ as above. Observe
that \begin{equation*}[(g,0),(\id,v)]=(\id,v\circ
g-v)\end{equation*} for $g\in\H(\rz^m)$ and $v\in\C_c(\rz^m)$.
Therefore
\begin{align*} &(\id,u)=(\id,\bar
w-w)=(\id,\bar w)(\id,-w)\\
&=(\id,v_1^--v^-_1\ci\phi^-+v_1^+-v^+_1\ci\phi^+)
(\id,v^-_2\ci\psi^--v_2^-+v^+_2\ci\psi^+-v_2^+)\\ &=
(\id,v_1^-\ci\phi^--v_1^-)^{-1}(\id,v_1^+\ci\phi^+-v_1^+)^{-1}(\id,v_2^-\ci\psi^--v_2^-)
(\id,v_2^+\ci\psi^+-v_2^+)\\
&=[(\phi^-,0),(\id,v_1^-)]^{-1}[(\phi^+,0),(\id,v_1^+)]^{-1}[(\psi^-,0),(\id,v_2^-)]
[(\psi^+,0),(\id,v_2^+)].
\end{align*}
 Consequently, $(\id,u)$ is in
the commutator subgroup as well.

Thus we have shown

\begin {lem} The semi-direct product group
$\H(\rz^m)\t_{\tau}\C_c(\rz^m)$ is a perfect group.
\end{lem}

\begin{rem}  Observe that $\phi^*$ and $\psi^*$ need
not be $C^1$-diffeomorphisms even if so are $\phi^*_{n}$ and
$\psi^*_{n}$.  Consequently, the construction is not valid in the
$C^1$ category.
\end{rem}

\begin{rem} Analogous results to Lemma 2.1 were also clue ingredients
of proofs in [4], [1], and [2].  However the proofs of these
analogs were much easier than the present one as they made use of
a stability property for diffeomorphisms or Lipschitz
homeomorphisms. This property is no longer true in the topological
category.

\end{rem}

\section{The proof of Theorem 1.1}

 Let $P:\H_G(M)\r\H(\bar M)$ be a
homomorphism given by $P(h)(\bar x)=\pi(h(x))$, where $x\in M$ is
any element such that $\pi(x)=\bar x$. Let $h\in\H_G(M)$. In view
of [8] and [6] it suffices to consider $h\in\ker P$. In fact,
$P(h)$ can be decomposed as $P(h)=g_1\circ\cdots\circ g_r$ such
that $g_i\in\H(\bar M)$ is supported in a ball $B_i$ with
$\pi|\pi^{-1}(B_i)$ trivial,  $i=1\ld r$. Then each $g_i$ can be
obviously lifted to $h_i\in\H_G(M)$, i.e. $P(h_i)=g_i$. Due to [8]
each $h_i$ can be written as a product of commutators and $h\circ
h_r^{-1}\circ\cdots\circ h_1^{-1}\in\ker P$.

 Assume now that
$h\in\ker P$. Then $h$ can be identified with a continuous
function $\bar h:\bar M\r G$ such that $h(x)=x\cdot\bar h(\bar
x)$, where $\pi(x)=\bar x$.

Let us identify $\L(G)$, the Lie algebra of $G$, with $\rz^q$,
$q=\dim(G)$, by means of a basis $(X_1\ld X_q)$  of $\L(G)$. Let
$\Phi:\rz^q\supset V\r U\s G$ be a chart given by $\Phi(t_1\ld t
_q)=(\exp\,t_1X_1)\ldots(\exp\,t_qX_q)$. Suppose that $h$ is so
small that the image of $\bar h$ is in $U$. If we let $\tilde
h=\Phi^{-1}\circ\bar h=(u_1\ld u_q)$ then we may and do assume
(after an obvious fragmentation) that each $u_i$ is supported in
an open ball $B_i$ in $\rz^m$ and $\rz^m$ is identified with a
chart domain. We can extend the semi-direct product structure from
$\H(\rz^m)\t_{\tau}\C_c(\rz^m)$ to
$\H(\rz^m)\t_{\tau}\C_c(\rz^m,\rz^q)$, where $\C_c(\rz^m,\rz ^q)$
is the space of compactly supported $\rz^q$-valued functions, by
the formulae $(h,(v_1\ld v_q))=(\id,(v_1\ld v_q)\circ
h)\cdot(h,0)$ and $(\id,(v_1\ld v_q))=(\id,v_1)\cdots(\id,v_q)$.
In view of Lemma 2.1, each $(\id,u_i)$ is
 in the commutator subgroup of
$\H(\rz^m)\t_{\tau}\C_c(\rz^m,\rz^q)$. Consequently, $h$ is in the
commutator subgroup of $\H_G(M)$. This completes the proof.

\begin{rem} It is apparent that {\it Fragmentation
Property} holds for $\H_G(M)$: If $h\in\H_G(M)$ with $\supp(h)\s
B_1\cup\ldots\cup B_r$ then $h$ can be written as
$h=h_1\circ\cdots \circ h_s$, where $\supp(h_i)\s B_{j(i)}$,
$i=1\ld s$.
\end{rem}

\begin{cor} Let $M$ be a topological $G$-manifold with one orbit
type. Then $\H_G(M)$ is a perfect group.
\end{cor}

Indeed, if $H$ is the isotropy group of a point of $M$ then
$M^H=\{x\in M: \, H\,\,\hbox{fixes}\,\,x\}$ is a free
$N(H)/H$-manifold, where $N(H)$ is the normalizer of $H$ in $G$.
Since $\H_G(M)$ is isomorphic to $\H_{N(H)/H}(M^H)$ (cf. [6]),
Corollary 3.2 follows from Theorem 1.1.

\def\wyr#1{\textit{#1}}

\end{document}